\documentstyle[12pt]{article}

\newcommand{\cD}{{\cal D}}

\newcommand{\cF}{{\cal F}}

\newcommand{\cH}{{\cal H}}

\newcommand{\cK}{{\cal K}}

\newcommand{\cO}{{\cal O}}

\newcommand{\cT}{{\cal T}}

\newcommand{\cW}{{\cal W}}

\newcommand{\B}[1]{\begin{#1}}
\newcommand{\E}[1]{\end{#1}}

\newcommand{\aX}{\mbox{$X_{an}$}}

\newcommand{\C}{\mbox{$\bf C$}}     
\newcommand{\Z}{\mbox{$\bf Z$}}     
\newcommand{\Q}{\mbox{$\bf Q$}}

\newcommand{\im}{{\rm im}\,}

\newcommand{\df}{\mbox{\,$\stackrel{\pp{\rm def}}{=}$}\,}

\newcommand{\by}[1]{\stackrel{#1}{\rightarrow}}

\newcommand{\implies}{\mbox{$\Rightarrow$}}

\newcommand{\into}{\hookrightarrow}

\newcommand{\ie}{{\it i.e.\/}\ }
\newcommand{\eg}{{\it e.g.\/}\ }
\newcommand{\cf}{{\it cf.\/}\ }
\newcommand{\op}{{\it op.\/ cit.\/}\ }
\newcommand{\sZ}{\mbox{\scriptsize{$\Z$}}}

\newcommand{\pp}[1]{\mbox{$\scriptscriptstyle {#1}$}}

\newcommand{\limdir}[1]{{\displaystyle{\mathop{\rm
lim}_{\buildrel\longrightarrow\over{#1}}}}\,}

\title{On the Deligne--Beilinson cohomology
sheaves}

\author{by Luca {\sc Barbieri-Viale}}

\date{}

\begin{document}

\maketitle

\B{abstract}
We are showing that the Deligne--Beilinson cohomology
sheaves $\cH^{q+1}(\Z(q)_{\cD})$ are torsion free by
assuming Kato's conjectures hold true for function fields.
This result is `effective' for $q=2$; in this case,
by dealing with `arithmetic properties' of the presheaves
of mixed Hodge structures defined by singular cohomology,
we are able to give a cohomological characterization of the
Albanese kernel for surfaces with $p_g=0$.
\E{abstract}

\tableofcontents

\section*{Introduction} For $X$ a compact complex
algebraic manifold the Deligne cohomology
$H^*(X,\Z(\cdot)_{\cD})$ is defined by taking the
hypercohomology of the truncated De Rham complex
augmented over $\Z$. The very extension of such a
cohomology theory to arbitrary algebraic complex
varieties is usually called Deligne--Beilinson cohomology
(\eg see \cite{GI} for definitions and properties
or \cite{HV} for details). The associated Zariski sheaves
$\cH^*(\Z(\cdot)_{\cD})$ have groups of global sections
which are birational invariants of smooth complete
varieties (see \cite{BV1}, \cite{BV2}). The motivation
for this paper is to investigate these invariants. We can
show that the Deligne--Beilinson cohomology sheaves
$\cH^{q+1}(\Z(q)_{\cD})$ are {\it torsion
free\,} by assuming Kato's conjectures hold true for
function fields (see \S 2). In particular
$\cH^{3}(\Z(2)_{\cD})$ is actually torsion free thanks to
Merkur'ev--Suslin's Theorem on $K_2$. Thus these
invariants {\it vanish\,} for unirational varieties. If
only $H^2(X,\cO_X)=0$ we then can show (see \S 3) that the
group of global sections of $\cH^3(\Z(2)_{\cD})$ is {\it
exactly\,} the kernel of the cycle map $CH^2(X)\to
H^4(X,\Z(2)_{\cD})$ in Deligne cohomology \ie the kernel
of the Abel--Jacobi map.  This fact generalizes the result
of H.Esnault for $0$-cycles (\cf \cite[Theorem~2.5]{H}),
in the case of codimension $2$ cycles, to $X$ of arbitrary
dimension and it is obtained by a different proof.
Concerning the discrete part $\cF^{2,2}_{\sZ}$ of the
Deligne--Beilinson cohomology sheaf $\cH^2(\Z(2)_{\cD})$
(as defined in \op or \S 2) we can
describe, for any $X$ proper and smooth, the torsion of
$H^1(X,\cF^{2,2}_{\sZ})$ in terms of `trascendental
cycles' and $H^3(X,\Z)_{tors}$, and we can see that there
are not non-zero global sections of it (see \S 3).  For
surfaces with $p_g=0$ we are able to compute the group of
global sections of $\cH^3(\Z(2)_{\cD})$ ---  Bloch's
conjecture is that $H^0(X,\cH^3(\Z(2)_{\cD}))=0$ --- by
means of some short exact sequences (see \S 4) involving
the discrete part $\cF^{2,2}_{\sZ}$ and the Hodge
filtration. In order to do that we are firstly arguing
(see \S 1) with arithmetic resolutions of the Zariski
sheaves associated with the presheaves of mixed Hodge
structures defined by singular cohomology: the Hodge and
weight filtrations do have corresponding coniveau spectral
sequences, the $E_2$ terms of which are given by the
cohomology groups of the Zariski sheaves associated to
such filtrations.\\ I would like to thank B.Kahn for
friendly useful conversations on some of the matters
contained herein.

\section*{Notations}

Throughout this note $X$ is a complex algebraic variety.
We will denote by $H^*(X,A)$ (resp. $H_*(X,A)$) the
singular cohomology (resp. Borel--Moore homology) of the
associated analytic space $\aX$ with coefficients in
$A$ where $A$ would be $\Z,\Z/n, \C$ and $\C^*$; we
let $H^*(X)$ (resp. $H_*(X)$)  be the corresponding mixed
Hodge structure (see \cite{D}). We will denote by
$W_iH^*(X)$ (resp. $W_{-i}H_*(X)$) the $\Q$-vector spaces
given by the weight filtration and by $F^iH^*(X)$ (resp.
$F^{-i}H_*(X)$) the real vector spaces given by the Hodge
filtration. For the ring $\Z$ of integers we will denote
by $\Z(r)$ the Tate twist in Hodge theory and by
$H^*(X,\Z(r)_{\cD})$ the Deligne--Beilinson cohomology
groups (see \cite{HV}, \cite{GI}). The Tate twist induces
the twist $A\otimes\Z(r)$ in the coefficients which we
will denote $A(r)$ for short. We will denote by
$\cH^*(A(r))$ and $\cH^*(\Z(r)_{\cD})$ the Zariski sheaves
on a given $X$ associated to singular cohomology and
Deligne--Beilinson cohomology respectively.

\section{Arithmetic resolutions in mixed Hodge theory}

Let $Z\into X$ be a closed subscheme of the complex
algebraic variety $X$. According with Deligne \cite[8.2.2
and 8.3.8]{D} the relative cohomology groups $H^*_Z(X,\Z)$
(= $H^*(X {\rm mod} X-Z,\Z)$ in \op) carry out a mixed
Hodge structure fitting into long exact sequences
 \B{equation}\label{loc}
\cdots \to H_Z^j(X)\to H_T^j(X)\to
H_{T-Z}^j(X-Z)\to H_Z^{j+1}(X)\to \cdots
\E{equation}
for any pair $Z\subset T$ of closed subschemes of $X$.
As it has been remarked in \cite{JA} the assignation
$$Z\subseteq X \leadsto (H_Z^*(X),H_*(Z))$$
yields a Poincar\'e duality theory with
supports (see \cite{BO} and furthermore we have that
the above theory is appropriate for albegraic cycles in
the sense of \cite{BV2}) with values in the abelian
tensor category of mixed Hodge structures. In particular,
by considering the presheaf of vector spaces $$U\leadsto
F^iH^j(U)$$  (resp. $U\leadsto W_iH^j(U)$) and sheafifying
it on a fixed variety $X$, we obtain Zariski sheaves
$\cF^i\cH^j$ (resp. $\cW_i\cH^j$) filtering the sheaves
$\cH^j(\C)$. We then have: \B{prop}\label{arifilt}
Let $X$ be smooth. The `arithmetic resolution'
$$0\to\cH^q(\C)\to \coprod_{x\in X^0}^{} i_x H^{q}(x) \to
\coprod_{x\in X^1}^{} i_x H^{q-1}(x) \to \cdots\to
\coprod_{x\in X^q}^{} i_x\C\to 0$$
is a bifiltered quasi-isomorphism
$$(\cH^q(\C),\cF,\cW)\by{\simeq}(\coprod_{x\in
X^{\odot}}^{} i_x H^{q-\odot}(x),\coprod_{x\in
X^{\odot}}^{} i_x F, \coprod_{x\in X^{\odot}}^{}
i_x W)$$ yielding flasque resolutions
$$0\to gr^i_{\cF}gr_j^{\cW}\cH^q(\C)\to\coprod_{x\in
X^{0}}^{} i_x gr^i_{F}gr_j^{W}H^{q}(x)\to\cdots\to
\coprod_{x\in
X^{q}}^{} i_x gr^{i-q}_{F}gr_{j-2q}^{W}H^{0}(x)\to 0$$
\E{prop} \B{proof} Because of \cite[Theor.1.2.10
and 2.3.5]{D} the functors $F^n$, $W_n$ and $gr^n_F$ (any
$n\in\Z$) from the category of mixed Hodge structures to
that of real vector spaces are exact; $gr_n^W$ is exact as
a functor from mixed Hodge structures to pure $\Q$-Hodge
structures.  So the claimed results are
obtained via the `locally homologically effaceable'
property (see \cite[Claim p. 191]{BO}) by construction of
the arithmetic resolution (granted by
\cite[Theor.4.2]{BO}). For example: by applying $F^i$ to
the long exact sequences (\ref{loc}), taking direct limits
over pairs $Z\subset T$ filtered by codimension and
sheafifying, we do obtain a flasque resolution of length
$q$ $$0\to \cF^i\cH^q\to \coprod_{x\in X^0}^{} i_x
F^iH^{q}(x) \to \coprod_{x\in X^1}^{} i_x
F^{i-1}H^{q-1}(x) \to \cdots$$ where $$F^*H^{*}(x)\df
\limdir{\pp{ U \mbox{\ open\ } \subset \overline{\{ x\}
}}} F^*H^*(U)$$ By this method we obtain as well a
resolution of $\cW_j$ $$0\to \cW_j\cH^q\to \coprod_{x\in
X^0}^{} i_x W_jH^{q}(x) \to \coprod_{x\in X^1}^{} i_x
W_{j-2}H^{q-1}(x) \to \cdots$$
These resolutions grant us of the claimed bifiltered
quasi-isomorphism.  (Note: for $X$ of dimension
$d$, the fundamental class  $\eta_X$ belongs to
$W_{-2d}H_{2d}(X)\cap F^{-d}H_{2d}(X)$ so that `local
purity' yields the shift by two for the weight filtration
and the shift by one for the Hodge filtration). In the same
way we obtain resolutions of $gr^i_{\cF}$, $gr_j^{\cW}$
and $gr^i_{\cF}gr_j^{\cW}$.
\E{proof}
We may consider the twisted Poincar\'e duality
theory $(F^nH^{*},F^{-m}H_*)$ where the integers $n$ and
$m$ play the role of twisting  and indeed we have
$$F^{d-n}H^{2d-k}_Z(X)\cong F^{-n}H_k(Z)$$ for $X$ smooth
of dimension $d$. Via the arithmetic resolution of
$\cF^i\cH^q$ we then have the following: \B{cor} Let
assume $X$ smooth and let $i$ be a fixed integer. We then
have a `coniveau spectral sequence'
\B{equation}\label{conifilt} E^{p,q}_2 =H^p(X,\cF^i\cH^q)
\implies F^iH^{p+q}(X) \E{equation} where
$H^p(X,\cF^i\cH^q)=0$ if $q<\mbox{min}(i,p)$.
\E{cor}
\B{rmk} Concerning the Zariski sheaves
$gr^i_{\cF}\cH^q$ and $\cH^q/\cF^i$ we
indeed obtain corresponding coniveau spectral sequences
as above.
\E{rmk}
Because of the maps of `Poincar\'e
duality theories' $F^iH^*(-)\to H^*(-,\C)$ we do have as
well, maps of coniveau spectral sequences; on the
$E_2$-terms the map
$$H^p(X,\cF^i\cH^q)\to
H^p(X,\cH^q(\C))$$ is given by taking Zariski cohomology of
$\cF^i\cH^q \into \cH^q(\C)$. For example: if $i<p$ we
clearly have (by comparing the arithmetic resolutions):
$H^p(X,\cF^i\cH^p)\cong H^p(X,\cH^p(\C))$ and
$$H^p(X,\cH^p(\C))\cong NS^p(X)\otimes\C$$ by
\cite[7.6]{BO} where $NS^p(X)$ is the group of cycles of
codimension $p$ modulo algebraic equivalence. For $i=p$ we
still have: \B{thm}\label{Nero} Let $X$ be smooth. Then
$$H^p(X,\cF^p\cH^p)\cong NS^p(X)\otimes\C$$ \E{thm}
\B{proof} Because of the Proposition~\ref{arifilt}
$$H^p(X,\cF^p\cH^p)\cong {\rm coker} (\coprod_{x\in
X^{p-1}}^{}F^1H^{1}(x) \to \coprod_{x\in X^p}^{} \C)$$
whence the canonical map $H^p(X,\cF^p\cH^p)\to
NS^p(X)\otimes\C$ is surjective. To show the injectivity,
via the arithmetic resolution we see that
$$H^{2p-1}_{Z^{p-1}}(X,\C)\cong\mbox{ker}(\coprod_{x\in
X^{p-1}}^{}H^{1}(x) \to \coprod_{x\in X^p}^{} \C)$$
where $H^{*}_{Z^{i}}$ denotes the direct limit of the
cohomology groups with support on closed subsets of
codimension $\geq i$; indeed this formula is obtained
by taking the direct limit of (\ref{loc}) over pairs
$Z\subset T$ of codimension $\geq p$ and $\geq p-1$
respectively, since  $H^{2p-1}_{Z^{p}}=0$ and
$$H^{2p}_{Z^{p}}(X,\C) = \coprod_{x\in X^p}^{} \C$$
Furthermore
$$F^pH^{2p-1}_{Z^{p-1}}\cong\mbox{ker}(\coprod_{x\in
X^{p-1}}^{}F^1H^{1}(x) \to \coprod_{x\in X^p}^{} \C)$$
and
$$H^{2p-1}_{Z^{p-1}}/F^p\cong \coprod_{x\in
X^{p-1}}^{}gr_F^0H^{1}(x)$$
since the arithmetic resolution of $\cH^p/\cF^p$ has
lenght $p-1$. Thus we have that
$$\mbox{image}(\coprod_{x\in
X^{p-1}}^{}F^1H^{1}(x) \to \coprod_{x\in X^p}^{} \C)=
\mbox{image}(\coprod_{x\in
X^{p-1}}^{}H^{1}(x) \to \coprod_{x\in X^p}^{} \C)$$
 \E{proof}
\B{rmk} Note that by considering the sheaf
$\cH^q(\C)$(=$\cF^0\cH^q$) on $X$ filtered by the
subsheaves $\cF^i\cH^q$ we have as usual (\cf
\cite[1.4.5]{D}) a spectral sequence
$${}_{\cF}E_1^{r,s}=H^{r+s}(X,gr^s_{\cF}\cH^q)\implies
H^{r+s}(X,\cH^q(\C))$$ with induced `aboutissement'
filtration $$F^iH^p(X,\cH^q) \df \im (H^p(X,\cF^i\cH^q)\to
H^p(X,\cH^q(\C)))$$ By the above Theorem we can see that
$F^iH^p(X,\cH^p)$ is uninteresting, since it gives the all
N\'eron--Severi group if $i\leq p$ and vanishes otherwise.
\E{rmk}
\B{rmk} As an immediate consequence of this Theorem, via
the coniveau spectral sequence (\ref{conifilt}), we see the
well known fact that the cycle map
$c\ell^p:NS^p(X)\otimes\C\to H^{2p}(X,\C)$ has its image
contained in $F^pH^{2p}(X)$.
\E{rmk}
For any $X$ smooth and proper, we have that
$F^2H^2(X)=H^0(X,\cF^2\cH^2)=H^0(X,\Omega^2_X)$ and
\B{equation}\label{genus}
 H^0(X,\cH^2/\cF^2)\cong
H^0(X,\cH^2(\C))/H^0(X,\Omega^2_X)\cong
\frac{H^2(X,\C)}{H^0(X,\Omega^2_X)\oplus NS(X)\otimes\C}
\E{equation} where $H^0(X,\cH^2(\C))\cong
H^0(X,\cH^2(\Z))\otimes\C$ and $H^0(X,\cH^2(\Z)) = \im
(H^2(X,\Z)\to H^2(X,\cO_X))$ is the lattice of
`trascendental cycles'.  The formula (\ref{genus}) can be
obtained, for example, by the exact sequence
(given by the coniveau spectral sequence since $\cF^2\cH^1
=0$) $$0\to H^1(X,\cH^1(\C))\to H^2(X)/F^2\to
H^0(X,\cH^2/\cF^2)\to 0$$ because of
$H^1(X,\cH^1(\C))=NS(X)\otimes\C$.

\section{Deligne--Beilinson cohomology sheaves}

Let $X$ be smooth over $\C$. Let consider the Zariski
sheaf $\cH^*(\Z(r)_{\cD})$ associated to the presheaf of
Deligne--Beilinson cohomology groups $U\leadsto
H^*(U,\Z(r)_{\cD})$ on $X$. We have canonical long exact
sequences of sheaves on $X$
\B{equation}\label{modf} \cdots\to
\cH^q(\Z(r))\to \cH^q(\C)/\cF^r \to
\cH^{q+1}(\Z(r)_{\cD}) \to \cH^{q+1}(\Z(r))\to\cdots
\E{equation} \B{equation}\label{plusf} \cdots\to
\cH^{q}(\Z(r)_{\cD}) \to \cH^{q}(\Z(r))\oplus
\cF^r\cH^q\to \cH^q(\C)\to \cH^{q+1}(\Z(r)_{\cD})\to\cdots
\E{equation}
\B{equation}\label{star}\cdots\to
\cF^r\cH^{q}\to\cH^q(\C^*(r))\to
\cH^{q+1}(\Z(r)_{\cD}) \to \cF^r\cH^{q+1}\to\cdots
\E{equation}
obtained by sheafifying the usual long exact sequences
coming with Deligne--Beilinson cohomology (see
\cite[Cor.2.10]{HV}).
For example, if $r=0$ then
$\cH^q(\C)/\cF^0=0$ and (\ref{modf}) yields the
isomorphism $\cH^{*}(\Z(0)_{\cD}) \cong \cH^{*}(\Z)$;
(\ref{plusf}) splits in trivial short exact
sequences and (\ref{star}) give us the following short
exact sequence \B{equation}\label{universal}
0\to \cH^q(\Z)\to\cH^q(\C)\to\cH^q(\C^*)\to 0
\E{equation}
whenever the sheaves $\cH^q(\Z)$ and $\cH^{q+1}(\Z)$ are
torsion free which is the case if $q\leq 2$ (\cf
\cite[\S 3-4]{BV1} and \cite[p.1240]{BlS}). Torsion
freeness of $\cH^q(\Z)$ for all $q\geq 4$ is a conjectural
property (see \cite[\S 7]{BV1}). In order to show that
$\cH^{q+1}(\Z)$ is torsion free it sufficies to see that
$\cH^{q}(\Z)\to\cH^{q}(\Z/n)$ is an epimorphism for any
$n\in\Z$; via the canonical map  $\cO^*_X\to\cH^1(\Z)$ and
cup--product we obtain a map $(\cO^*_X)^{\otimes
q}\to\cH^q(\Z)$. The composition of $$(\cO^*_X)^{\otimes
q}\to\cH^q(\Z)\to \cH^{q}(\Z/n)$$ can be obtained as well
as (\cf \cite[p.1240]{BlS}) the composition of
$$(\cO^*_X)^{\otimes q}\by{sym} \cK_q^M \to \cH^q(\Z/n)$$
where by definition of Milnor's $K$-theory sheaf the
symbol map  $sym$ is an epimorphism thus we are left to
show that the Galois symbol $\cK_q^M \to \cH^q(\Z/n)$ is
an epimorphism (for the sake of exposition we are tacitly
fixing an $n$-th root of unity yielding a non-canonical
isomorphism $\cH^{q}_{\acute{e}t}(\mu_n^{\otimes q})\cong
\cH^q(\Z/n)$); this last map can be obtained by mapping
the Gersten complex for Milnor's $K$-theory to the
Bloch--Ogus arithmetic resolution of the sheaf
$\cH^q(\Z/n)$ \ie there is a commutative diagram
$$\B{array}{ccc} \cK_q^M &\to {\displaystyle
\coprod_{\eta\in X^0}^{}  i_{\eta} K^M_q(k(\eta))} &\to
{\displaystyle\coprod_{x\in X^1}^{} i_x
K_{q-1}^M(k(x))}\\ \downarrow & \downarrow &
\downarrow\\\cH^q(\Z/n) & \into
{\displaystyle\coprod_{\eta\in X^0}^{}
i_x H^{q}(\eta)}& \to
{\displaystyle\coprod_{x\in X^1}^{} i_x
H^{q-1}(x)}
\E{array}$$ where $H^*(\mbox{point})$ is the
Galois cohomology of  $k(\mbox{point})$. Indeed O.Gabber
announced the (universal) exactness of the above complex
of Milnor's $K$-groups. Thus: by assuming  Kato's
conjecture \ie $K^M_*(k(\mbox{point}))/n\cong
H^*(\mbox{point})$, a diagram chase yields the desired
projection $\cK_q^M \to \cH^q(\Z/n)$. (My thanks to B.Kahn
for having directed my attention to Gabber's result.) Let
consider the `discrete part' of the Deligne--Beilinson
cohomology sheaves which is by definition (\cf \cite[\S
1]{H}) $$\cF^{r,q}_{\sZ}\df \mbox{\ image\
}(\cH^q(\Z(r)_{\cD}) \to \cH^q(\Z(r)))$$ or, equivalently
by (\ref{modf}), the integral part of $\cF^r\cH^q$. We may
define the `trascendental part' of the Deligne--Beilinson
cohomology sheaves as follows $$\cT^{r,q}_{\cD}\df \mbox{\
kernel\ }(\cH^q(\Z(r)_{\cD}) \to \cH^q(\Z(r)))$$  In
particular, if $r=q$ we then have (\cf \cite[(1.3)$\alpha
)$]{H}) the short exact sequence (by (\ref{star}) or
(\ref{modf}) taking account of (\ref{universal}))
\B{equation} \label{disc} 0\to \cH^{q-1}(\C^*(q))\to
\cH^q(\Z(q)_{\cD})\to \cF^{q,q}_{\sZ}\to 0 \E{equation}
and moreover we have the following commutative diagram
with exact rows and columns
\B{equation}\B{array}{ccccccccc}\label{dia} &&0&&0&&0&&\\
&&\uparrow&&\uparrow&&\uparrow&&\\ 0&\to&
\cH^q(\Z(q))/\cF^{q,q}_{\sZ}&\to &\cH^q(\C)/\cF^q&\to &
\cT^{q,q+1}_{\cD}&\to &0\\
&&\uparrow&&\uparrow&&\uparrow&&\\
0&\to&\cH^q(\Z(q))&\to&\cH^q(\C)&\to&\cH^q(\C^*(q))
&\to&0\\ &&\uparrow&&\uparrow&&\uparrow&&\\
0&\to&\cF^{q,q}_{\sZ}&\to&\cF^q&\to&
\cF^q/\cF^{q,q}_{\sZ}&\to&0\\
&&\uparrow&&\uparrow&&\uparrow&&\\
&&0&&0&&0&&
\E{array}\E{equation}
where the middle row is given by (\ref{universal}), the
top one by (\ref{modf}) and the right-most column is
obtained by (\ref{star}).\B{lemma} We have a
short exact sequence of sheaves:
\B{equation}\label{unidel}
0\to\cH^q(\Z(r)_{\cD})/n\to \cH^q(\Z/n(r))\to
\cH^{q+1}((\Z(r)_{\cD})_{n-tors}\to 0
\E{equation} for all $q, r\geq 0$
and $n\in\Z$.\E{lemma}\B{proof} The sequence (\ref{unidel})
is obtained from the long exact sequence (\ref{modf}) as
follows. Since the sheaf $\cH^{q+1}(\Z(r))$ is torsion
free we do have that $$\cT^{r,q+1}_{\cD,n-tors}=
\cH^{q+1}((\Z(r)_{\cD})_{n-tors}$$
Since the sheaf $\cH^q(\C)/\cF^r$ is uniquely divisible we
have that  $$\cH^{q+1}((\Z(r)_{\cD})_{n-tors}=
(\cH^q(\Z(q))/\cF^{r,q}_{\sZ})\otimes \Z/n$$ because of
the following short exact sequence
$$0\to \cH^q(\Z(r))/\cF^{r,q}_{\sZ}\to \cH^q(\C)/\cF^r\to
\cT^{r,q+1}_{\cD}\to 0$$
Thus we get a short exact sequence
$$0\to \cF^{r,q}_{\sZ}/n \to \cH^q(\Z/n(r))\to
\cH^{q+1}((\Z(r)_{\cD})_{n-tors}\to 0$$
by tensoring with $\Z/n$ the canonical one induced by the
subsheaf $\cF^{r,q}_{\sZ}\into \cH^q(\Z(r))$. Since
$\cT^{r,q}_{\cD}$ is divisible we are done.\E{proof}
By considering the Bloch--Beilinson regulators
$$\cK_q^M\to \cH^q(\Z(q)_{\cD})$$ (simply obtained by the
fact that $\cK_1^M=\cO_X^*\cong \cH^1(\Z(1)_{\cD})$ and
cup--product) we have that the composition of
$$\cK_q^M\to \cH^q(\Z(q)_{\cD})\to \cH^q(\Z(q))\to
\cH^q(\Z/n(q))$$
is the Galois symbol (\cf \cite[\S 0 p.375]{HM}). We thus
have that the composition of
$$\cK_q^M/n\to\cH^q(\Z(q)_{\cD})/n\into\cH^q(\Z/n(q))$$
is an epimorphism if Kato's conjectures hold. Therefore,
by comparing with (\ref{unidel}), we have a proof of the
following: \B{thm}\label{Kato} Let assume that Kato's
conjectures (for function fields) hold true. On a smooth
$X$ we then have: \B{description}\item[{\it i)}] the sheaf
$\cH^q(\Z)$ is torsion free; \item[{\it ii)}] the sheaf
$\cH^{q+1}(\Z(q)_{\cD})$ is torsion free; \item[{\it
iii)}] there is a canonical isomorphism
$\cH^q(\Z(q)_{\cD})\otimes \Z/n\cong \cH^q(\Z/n(q))$
 \E{description}
for any $q\geq 0$.
\E{thm}
Note that the Theorem of Merkur'ev--Suslin (=
Kato's conjecture for $K^M_2$) ensures the previous
results when $q=2$. Moreover, by a standard argument (\cf
\cite[\S 2]{BV1} and \cite{BV2}) we have the following:
\B{cor} With the assumptions in the Theorem above, let
suppose that $X$ is moreover unirational and complete.
Then $$H^0(X,\cH^{q+1}(\Z(q)_{\cD}))=0$$
\E{cor}

\section{Coniveau versus Hodge filtrations}

We recall (see \cite{GI}) the
existence of arithmetic resolutions of the sheaves
$\cH^{*}(\Z(\cdot)_{\cD})$ thus the coniveau spectral
sequence \B{equation}\label{conideli}
{}_{\cD}E^{p,q}_2 =H^p(X,\cH^q(\Z(\cdot)_{\cD})) \implies
H^{p+q}(X,\Z(\cdot)_{\cD})
\E{equation}
and the formula (see
\cite{GI}) $H^p(X,\cH^p(\Z(p)_{\cD}))\cong CH^p(X)$. By the
spectral sequence (\ref{conideli}) we have a long exact
sequence  \B{equation}\label{delseq} 0\to
H^1(X,\cH^2(\Z(2)_{\cD}))\to H^3(X,\Z(2)_{\cD})\by{\rho}
H^0(X,\cH^3(\Z(2)_{\cD}))\by{\delta}CH^2(X) \E{equation}
The mapping $\delta$ is just a differential between
${}_{\cD}E_2$-terms of the coniveau spectral sequence
(\ref{conideli}); we still have  $$\mbox{image\, }\delta =
\ker (CH^2(X)\by{c\ell}H^4(X,\Z(2)_{\cD}))$$
\B{prop}\label{alb} Let $X$ be proper and smooth. Then
$$H^0(X,\cF^{2,2}_{\sZ})=0$$ the group
$H^1(X,\cF^{2,2}_{\sZ})$ is infinitely divisible and
$$H^1(X,\cF^{2,2}_{\sZ})_{tors}\cong
H^0(X,\cH^2(\Q/\Z(2)))$$ If $H^2(X,\cO_X)=0$ then
$$H^0(X,\cH^3(\Z(2)_{\cD})) \cong  \ker
(CH^2(X)\by{c\ell}H^4(X,\Z(2)_{\cD}))$$
\ie $\rho =0$ in (\ref{delseq}), and
$H^1(X,\cF^{2,2}_{\sZ})\cong H^3(X,\Z)_{tors}$.
\E{prop}
\B{proof} By the canonical map of `Poincar\'e duality
theories' $$H^{\sharp -1}(-,\C^*(\cdot))\to
H^{\sharp}(-,\Z(\cdot)_{\cD})$$
we do obtain a map of coniveau spectral sequences and the
following commutative diagram
\B{equation}\label{cocco}\B{array}{ccc}
0\to H^1(X,\cH^2(\Z(2)_{\cD}))&\to
H^3(X,\Z(2)_{\cD})&\by{\rho} H^0(X,\cH^3(\Z(2)_{\cD}))\\
\uparrow &\uparrow &\uparrow\\
0\to NS(X)\otimes\C^*(2)&\to H^2(X,\C^*(2))&\to
H^0(X,\cH^2(\C^*(2)))\to 0
\E{array}
\E{equation}
where $H^1(X,\cH^1(\C^*(2)))\cong NS(X)\otimes \C^*(2)$ and
the left-most map is induced by the short exact
sequence of sheaves (\ref{disc}) whence
$H^0(X,\cF^{2,2}_{\sZ})$ and $H^1(X,\cF^{2,2}_{\sZ})$
are respectively the kernel and the cokernel of
$NS(X)\otimes\C^*(2)\to H^1(X,\cH^2(\Z(2)_{\cD})))$
in fact: the cokernel is computed by the vanishing of
$H^2(X,\cH^1(\C^*(2)))$ and the kernel is obtained
because of $H^0(X,\cH^1(\C^*(2)))\cong H^1(X,\C^*(2))$,
$H^0(X,\cH^2(\Z(2)_{\cD}))\cong H^2(X,\Z(2)_{\cD})$ (note
that $\C^*(2)\cong\cH^1(\Z(2)_{\cD})$ is flasque) and
$H^1(X,\C^*(2))\cong H^2(X,\Z(2)_{\cD})$ since $X$ is
proper \ie $F^2H^2\into H^2(X,\C^*(2))$. Furthermore we
have that $H^3(X,\Z(2)_{\cD})\cong H^2(X,\C^*(2))/F^2H^2$
and $H^0(X,\cF^{2,2}_{\sZ})$ vanishes because $F^2H^2\cap
NS(X)\otimes\C^*(2)=0$. Since $H^0(X,\cH^3(\Z(2)_{\cD}))$
is torsion free (by the Theorem~\ref{Kato}) then $\im\rho$
is torsion free and infinitely divisible indeed, therefore
$$H^1(X,\cH^2(\Z(2)_{\cD}))\otimes \Q/\Z =
H^1(X,\cF^{2,2}_{\sZ})\otimes \Q/\Z =0$$
and by taking the torsion subgroups in the diagram
(\ref{cocco}) we obtain the assertion about the torsion
of $H^1(X,\cF^{2,2}_{\sZ})$.
In order to show the second part of the statement,
since $NS(X)\cong H^2(X,\Z)$, we then have (by the bottom
row of the diagram (\ref{cocco}) above) that
$H^0(X,\cH^2(\C^*(2)))\cong H^3(X,\Z)_{tors}$ and
its image in $H^0(X,\cH^3(\Z(2)_{\cD}))$ is equal to
the image of $\rho$, whence the image of $\rho$ is zero
since it is torsion free. Since $F^2H^2=0$ by diagram chase
we obtain the last claim. \E{proof}
\B{rmk} The group $H^0(X,\cH^2(\Q/\Z))$ is effectively
the extension of $H^0(X,\cH^2(\Z))\otimes \Q/\Z$ by
$H^3(X,\Z)_{tors}$ because $H^0(X,\cH^3(\Z))$ is
torsion free.
\E{rmk}
In order to detect elements in the misterious
group of global sections $H^0(X,\cH^3(\Z(2)_{\cD}))$ we
dispose of the image of $H^0(X,\cH^2/\cF^2)$, see
(\ref{genus}), which is the same (\cf the diagram
(\ref{dia})) as the image of $H^0(X,\cH^2(\Z))\otimes
\C/\Q(2)=H^0(X,\cH^2(\C^*(2)))\otimes \Q$. Unfortunately
these images cannot be the entire group, in general.
Indeed we have that whenever the map
$$H^0(X,\cH^2(\C^*(2)))\to H^0(X,\cH^3(\Z(2)_{\cD}))$$ is
surjective then $\rho$ is surjective in (\ref{delseq})
(because of (\ref{cocco})) whence the cycle map is
injective which is not the case in general (indeed for any
surface with $p_g\neq 0$ the cycle map is not injective by
Mumford \cite{MU}).

\section{Surfaces with $p_g=0$}

In the following we let $X$ denote a complex algebraic
surface which is smooth and complete. Let $A_0(X)$ be the
subgroup of cycles of degree zero. Let  $\phi: A_0(X)\to
J^2(X)$ be induced by the canonical mapping to the
Albanese variety. It is well known (see \cite[Theor.2 and
Cor.]{GI}) that $c\ell\mid_{A_0(X)} =\phi$, where $c\ell$
is the cycle map in Deligne cohomology. We have:
\B{equation}\label{cacca} H^1(X,\cH^3(\Z(2)_{\cD}))=0
\E{equation}   Indeed: the sheaf $\cH^4(\Z(2)_{\cD})$
vanishes on a surface (as it is easy to see via the exact
sequence (\ref{modf})) and by the spectral sequence
(\ref{conideli}) we have that $H^1(X,\cH^3(\Z(2)_{\cD}))$
is the cokernel of the cycle map
$CH^2(X)\by{c\ell}H^4(X,\Z(2)_{\cD})$ but
coker$(c\ell)=$coker$(A_0(X)\to J^2(X))=0$. Finally we
have that  $H^2(X,\cH^3(\Z(2)_{\cD}))\cong
H^5(X,\Z(2)_{\cD})=0$. Thus the only possibly non-zero
terms in the spectral sequence (\ref{conideli}) are: those
giving the  exact sequence (\ref{delseq}),
$H^0(X,\cH^2(\Z(2)_{\cD}))\cong H^2(X,\Z(2)_{\cD})$ and
$H^0(X,\cH^1(\Z(2)_{\cD}))=\C^*$. We know by Mumford (see
\cite{MU}) that $A_0(X)\cong J^2(X)$ \ie $\delta=0$  in
(\ref{delseq}), implies $p_g=0$ therefore
$H^0(X,\cH^3(\Z(2)_{\cD}))=0$ (by Proposition~\ref{alb}).
Conversely we have the following:\\[5pt] {\bf Bloch's
conjecture\ \ } {\it If $p_g=0$ then $A_0(X)\cong
J^2(X)$.}\\[4pt] In order to test this conjecture we are
reduced to compute the uniquely divisible group
$H^0(X,\cH^3(\Z(2)_{\cD}))$ (\cf \cite{H}, \cite{GI}). We
can show the following: \B{prop}\label{compute} Let $X$ be
a smooth complete surface with $p_g=0$. We then have the
following canonical short exact sequences
\B{equation}\label{a} 0\to H^0(X,\cH^3(\Z(2)_{\cD}))\to
H^1(X,\cF^2/\cF^{2,2}_{\sZ})\to H^3(X,\C^*(2))\to 0
\E{equation} \B{equation}\label{b} 0\to
H^0(X,\cH^3(\Z(2)_{\cD}))\to
H^1(X,\cH^2(\Z(2))/\cF^{2,2}_{\sZ}))\to H^3/F^2\to 0
\E{equation} where \B{equation}\label{c} 0\to F^2H^3\to
H^1(X,\cF^2/\cF^{2,2}_{\sZ})\to
 A_0(X)\to 0
\E{equation}
\B{equation}\label{d}
0\to H^3(X,\Z)/tors\to
H^1(X,\cH^2(\Z(2))/\cF^{2,2}_{\sZ})\to A_0(X)\to 0
\E{equation}
 \E{prop}
\B{proof} All these exact sequences are obtained by
considering the exact diagram of cohomology groups
associated with the diagram of sheaves (\ref{dia}) (where
$\cT^{2,3}_{\cD} = \cH^3(\Z(2)_{\cD})$ on a surface)
taking account of Theorem~\ref{Kato}, Theorem~\ref{Nero},
Proposition~\ref{alb} and the coniveau spectral sequence
(\ref{conifilt}). For example, the sequence (\ref{a}) is
obtained by taking the long exact sequence of cohomology
groups associated with the right-most column of
(\ref{dia}), the fact that $H^1(X,\cH^2(\C^*(2)))\cong
H^3(X,\C^*(2))$ on a surface and the formula
(\ref{cacca}). For (\ref{b}) one has to use the top row
of (\ref{dia}), the formulas (\ref{genus}) and
(\ref{cacca}), and the fact that $H^3/F^2\cong
H^1(X,\cH^2/\cF^2)$. The left-most column of (\ref{dia})
yields (\ref{d}) because of $H^2(X,\cF^{2,2}_{\sZ})\cong
H^2(X,\cH^2(\Z(2)_{\cD}))\cong CH^2(X)$ by (\ref{disc})
(\cf \cite[1.3]{H}) and the map of sheaves
$\cH^2(\Z(2)_{\cD})\to \cH^2(\Z(2))$ induces the degree
map on $H^2$. For (\ref{c}) one has to argue with the
commutative square in the left bottom corner of
(\ref{dia}) and the isomorphim $H^2(X,\cF^2)\cong
H^2(X,\cH^2(\C))\cong \C$: remember that
$H^1(X,\cF^{2,2}_{\sZ})=H^3(X,\Z)_{tors}$ whence it goes to
zero in $F^2H^3\cong H^1(X,\cF^2\cH^2)$.
\E{proof}
\B{rmk} Because of this Proposition,  Bloch's
conjecture is equivalent to showing that the
canonical injections of sheaves
$\cF^2/\cF^{2,2}_{\sZ}\into \cH^2(\C^*(2))$ or
$\cH^2(\Z(2))/\cF^{2,2}_{\sZ}\into  \cH^2/\cF^2$ remain
injections on $H^1$. It would be very nice to know of any
reasonable description of the Zariski cohomology classes
of these subsheaves.  \E{rmk}


\begin{thebibliography}{}

\bibitem{BV1}{\sc L.Barbieri-Viale}: Cicli di
codimensione $2$ su variet\`a unirazionali complesse, to
appear in the {\it Ast\'erisque} Volume of Proc. of
the $K$-theory Conf. in Strasbourg, 1992.

\bibitem{BV2}{\sc L.Barbieri-Viale}: $\cH$-cohomologies
versus algebraic cycles, preprint, 1994.

\bibitem{BO}{\sc S.Bloch~}and{\sc~A.Ogus}: Gersten's
conjecture and the homology of schemes, {\it Ann.Sci.Ecole
Norm.Sup.} {\bf 7} (1974) 181--202.

\bibitem{BlS}{\sc S.Bloch}~and{\sc~V.Srinivas}: Remarks on
correspondences and algebraic cycles, {\it Amer. J.
Math.} {\bf 105} (1983), 1235-1253.

\bibitem{D}{\sc P.Deligne}: Th\'eorie de Hodge
II,III , {\it Publ. Math.}\, IHES {\bf 40} (1972) 5--57
and {\bf 44} (1974) 5--78.

\bibitem{HM}{\sc H.Esnault}: Une remarque sur la
cohomologie du faisceau de Zariski de la $K$-th\'eorie de
Milnor sur une vari\'et\'e lisse complexe, {\it Math.Z.}
{\bf 205} (1990) 373--378.

\bibitem{H}{\sc H.Esnault}: A note on the cycle map, {\it
J.reine angew.Math.} {\bf 411} (1990) 51--65.

\bibitem{HV}{\sc H.Esnault}~and{\sc~E.Viehweg}:
Deligne-Beilinson cohomology, in  Beilinson's
conjectures on special values of $L$-functions, {\it
Academic Press Perspectives in Math.} {\bf 4}
(1988) 43--92.

\bibitem{GI}{\sc H.Gillet}: Deligne homology and the
Abel--Jacobi maps, {\it Bull.AMS} {\bf 10}
(1984) 285--288.

\bibitem{JA}{\sc U.Jannsen}: Mixed motives and algebraic
$K$-theory, {\it Springer} LNM {\bf 1400}, 1990.

\bibitem{MU}{\sc D.Mumford}: Rational equivalence of
$0$-cycles on surfaces, {\it J.Math.Kyoto Univ.} {\bf 9}
(1968) 195--204.

\end{thebibliography}
 \end{document}